\documentclass[11pt]{article}

\usepackage{graphics}
\usepackage{amsfonts}
\usepackage{amssymb}
\usepackage{amsmath}
\usepackage{amsthm}
\usepackage{url}
\usepackage{color}
\usepackage{multirow}
\usepackage{enumerate}
\usepackage{hyperref}

\usepackage{algorithm}
\usepackage{algorithmic}

\theoremstyle{definition}%an alternative is \theoremstyle{remark}

\theoremstyle{definition}

\theoremstyle{definition}

\newcommand{\R}{\mathbb{R}}

\title{On solving the densest $k$-subgraph problem on large graphs}
\author{Renata Sotirov\thanks{Department of Econometrics and OR, Tilburg University, The Netherlands. {\tt r.sotirov@uvt.nl} }
 }

\date{}

\begin{document}
\maketitle

\begin{abstract}
The densest $k$-subgraph problem  is the problem of finding a  $k$-vertex subgraph  of a graph  with the maximum number of edges.
In order to solve large instances of the  densest $k$-subgraph problem,
we introduce two algorithms that are based on the random coordinate descent approach.
Although it is common use to update  at most two random  coordinates  simultaneously in each iteration of an algorithm, our algorithms  may simultaneously update many coordinates.
We show the benefit of updating more than two coordinates  simultaneously  for solving the  densest $k$-subgraph problem, and
 solve large problem instances  with up to $2^{15}$ vertices.
\end{abstract}


\section{Introduction}\label{sec:intro}

The densest $k$-subgraph (D$k$S) problem is the problem of finding a subgraph  of the given graph with exactly $k$ vertices such that the number of edges in the subgraph is maximal.
The densest $k$-subgraph problem is known in the literature under various names, including
the heaviest unweighted subgraph problem \cite{kortsarz1993choosing}, the $k$-cluster problem \cite{corneil1984clustering}, or the $k$-cardinality subgraph problem \cite{bruglieri2006annotated}.
The densest $k$-subgraph problem can be seen as a special case of the maximum $k$-dispersion problem \cite{ravi1994heuristic}.
The maximum $k$-dispersion problem  is the problem of finding  $k$ vertices in a graph that
maximize a function of the distances between the chosen  vertices.
The D$k$S problem can also be seen as a special case of the heaviest $k$-subgraph problem, which
is the problem of finding a subgraph with $k$ vertices that maximizes the sum of the edge weights in the subgraph.

The D$k$S problem is known to be NP-hard. In \cite{feige1997densest} it is proven that the problem is NP-hard  for graphs whose  maximum degree is equal to three.
The densest $k$-subgraph problem  is NP-hard even for very restricted classes of graphs,
such as  bipartite and chordal graphs \cite{corneil1984clustering}, or planar graphs \cite{Keil}.
However, it is trivial on trees. %% because any subtree of $k$ vertices contains exactly $k-1$ edges.
The D$k$S problem  is solvable in polynomial time on graphs whose maximum degree  is equal to two,
as well as on cographs, split graphs, and $k$-trees, see \cite{corneil1984clustering}.

There are many applications of the problem.
The densest $k$-subgraph problem  plays a role in analyzing web graphs and different social networks.
Namely, one of the main challenges for web search engines is the detection of link spams, see Henzinger et al.~\cite{henzinger2003challenges}.
Link spams are websites that are linked to each other in order to manipulate the search engine rankings.
Many of the dense subgraphs in web graphs are link spams.
Gibson et al.~\cite{gibson2005discovering} propose an algorithm that extracts dense subgraphs  in huge graphs in order to identify link spams.
Angel et al.~\cite{Angel12} analyze social networks to identify  real-time stories by searching for dense subgraphs of the given size.
The heaviest $k$-subgraph problem can be used to create interest groups of people.
For example, when  organizing an opening party where  participants should be as similar as possible, or to analyze political vote data, see  \cite{tsourakakis2015k}.
The heaviest $k$-subgraph may be used to find teams of employees with the highest collaborative compatibility, see \cite{pragarauskas2012multi}.\\

 We next list problems that are related to the D$k$S problem.
The following two versions  of the densest $k$-subgraph problem are introduced in \cite{andersen2009finding}:
the densest at-least-$k$-subgraph problem and the densest at-most-$k$-subgraph problem.
The densest at-least-$k$-subgraph (resp.~the densest at-most-$k$-subgraph) problem is the problem of finding an induced subgraph of highest average degree with at least (resp.~at most) $k$ vertices.
Andersen  and  Chellapilla \cite{andersen2009finding}  present  an efficient $1/3$-approximation algorithm for  the densest at-least-$k$-subgraph problem.
There are no efficient approximation algorithms for the densest at-most-$k$-subgraph problem.
The problem of finding a subgraph of maximum node weights with exactly $k$ edges is considered in \cite{goldschmidt1997k}.
The sparsest $k$-subgraph problem finds the subgraph with $k$ vertices and the minimum number of edges, see  \cite{bougeret2014parameterized}.\\

{\bf Outline and main results.}  This paper is structured as follows.
Section \ref{solving approaches} presents an integer programming formulation of the problem and lists various solving approaches from the literature.
In Section \ref{methods overview} we provide an overview of recently introduced  methods for solving large scale optimization problems.
We present two new algorithms for solving  the D$k$S problem in Section \ref{sect:onSolveDkS}.
Our algorithms are applied to the relaxation of the  D$k$S problem,  see \eqref{lpQP}.
Our first algorithm considers a quadratic optimization subproblem with linear constraints, and the second one a linear programming optimization subproblem.
The main difference between the here presented algorithms and those in the literature is that we allow updating more than two random coordinates
 simultaneously in each iteration of our algorithms.
 We show here that for an appropriate number of simultaneous updates, our second algorithm converges to an integer solution vector (!).
 This convergence is not proven with theoretical convergence analysis, but only empirically observed.

Our extensive numerical results show  that we find densest subgraphs in large graphs in short time, see Section \ref{sect:numer}.
For example,  we find a densest subgraph with 25 vertices  in a graph with 23,133 vertices and 93,497 edges in less than 4 minutes.
 Since the densest subgraph is a clique in this case,  we know that we found an optimal subgraph.
Exact approaches for finding densest $k$-subgraphs  can not cope with graphs that have more than 160 vertices.
On the other hand, the best  heuristic approaches are tested on random instances with at most 3,000 vertices.
We test our algorithms on  real-world  data and on randomly generated data with up to 32,768 vertices.

\section{The problem formulation and {solution approaches}} \label{solving approaches}

The densest $k$-subgraph problem can be formulated as a  quadratic optimization problem with binary variables.
Let  $G=(V,E)$ be an undirected graph with vertex set $V$, $|V|=n$, and edge set $E$,  $|E|=m$.
 Let $k$ be a  positive integer between 3 and $n-2$,
and $x_{i}$ a binary variable that obtains value one if vertex $i$ is in the densest $k$-subgraph and zero otherwise.
We denote by $A$ the adjacency matrix of $G$.

The densest $k$-subgraph problem can be formulated as follows:
\begin{equation}\label{QP}
\begin{array}{rl}
\max & x^{\mathrm T}Ax \\[1ex]
{\rm s.t.} & \sum\limits_{i=1}^{n} x_{i} = k \\[1ex]
& x_{i} \in \{0,1\}, ~~ \forall i \in \{1,...,n\}.
\end{array}
\end{equation}
In the sequel, we list  approaches that are used for solving the D$k$S.
Billionnet \cite{billionnet2005different} derived four different mixed-integer linear programming formulations for the heaviest $k$-subgraph problem,
and three different  mixed-integer linear programming formulations  for  the densest $k$-subgraph problem.
Numerical results in \cite{billionnet2005different} show that the quality of a formulation is related to the density of an instance.
In \cite{billionnet2009improving}, the authors solve instances of the densest $k$-subgraph problem
by reformulating the  non-convex quadratic  problem  \eqref{QP} into an equivalent problem with a convex objective function.
Such reformulation requires solving an associated semidefinite programming problem.
The reformulated problem is then solved by using a branch-and-bound algorithm.
This approach is tested on random graphs with at most 100 vertices.
Numerical results show that the proposed convexification approach improves efficiency of the branch and bound algorithm.
However, solving the related  semidefinite program may be costly.

Malick and Roupin \cite{Malickroupin} solve instances of the D$k$S problem to optimality using semidefinite programming.
Namely, they solve  a semidefinite programming problem in each node of a branch and bound tree. The largest solved instances of the D$k$S has  120  vertices.
{Krislock, Malick and Roupin  \cite{KrisMalRoup:16} report solving  hard instances of the D$k$S problem  with up to 160 vertices   by using a
semidefinite branch and bound algorithm.}
Semidefinite programming relaxations are also used in the design of approximation algorithms for the D$k$S, see e.g., \cite{FeigeLeagberg01,{YeZhang}}.
{One can find an overview of SDP relaxations for the densest $k$-subgraph problem in \cite{Rendl15}.}
The above mentioned results show that it is extremely difficult to find a densest $k$-subgraph in a graph that has more than 160 vertices by using exact approaches.

{
A number of recent results have focused on  recovering  planted $k$-subgraphs by using convex relaxation techniques, see e.g., \cite{AmesVavas,Ames}.
Ames and Vavasis \cite{AmesVavas} show that the maximum clique in a graph consisting of a single large clique can be identified
 from the minimum nuclear norm solution of a particular system of linear inequalities.
Ames \cite{Ames}  establishes analogous recovery guarantees for a convex relaxation of the planted clique problem that is robust to noise.
For a survey on the topic see Li et al., \cite{LiChenXu}.
 }

In 2001, Feige at al.~\cite{feige2001dense} provide an approximation algorithm  for the D$k$S problem with approximation ratio of $n^{\delta -\epsilon}$ for some small $\epsilon$.
In \cite{bhaskara2010detecting}, it is presented an approximation algorithm that for every $\epsilon> 0$ approximates the D$k$S problem within a ratio of  $n^{1/4+\epsilon}$ in  $n^{O(1/\epsilon)}$ time.
The most recent results on the superpolynomial approximation algorithms for the  D$k$S one can find in   \cite{Bourgeois17}.
In  \cite{khot2006ruling}, Khot proves  that there does not exist a polynomial time approximation scheme (PTAS) for the densest $k$-subgraph problem in general graphs.
However, there exist   polynomial time approximation schemes for a few special problem cases.
Arora et al.~\cite{arora1995polynomial} provide a PTAS for the D$k$S problem on dense instances. Nonner \cite{nonner2016ptas} drives a PTAS  for interval graphs.

Different heuristic methods are tested for solving the densest $k$-subgraph problem.
 Kincaid \cite{kincaid1992good} uses  simulated annealing and tabu search heuristics to solve the D$k$S problem.
His results show that  the tabu search algorithm  performs better than the simulated annealing algorithm  for solving the densest $k$-subgraph  problem.
In  \cite{macam2002}, Macambira implements  tabu search heuristics  for the heaviest $k$-subgraph problem.
Although the tabu search algorithm from  \cite{macam2002} does not perform diversification, it  outperforms the greedy randomized adaptive search procedure.
A variable neighborhood search (VNS) heuristics  for the heaviest subgraph problem and graphs up to  3,000 vertices is implemented by  Brimberg et al.~\cite{brimberg2009variable}.
Their results show that the VNS outperforms the tabu search heuristic and multi-start local search heuristics in solving the D$k$S.
The VNS performs extremely well on sparse graphs.
Running times needed to find the best solutions for instances with 3,000 vertices is about 425 seconds.
A heuristic based on a two-step filtering approach is used to extract dense web communities in Dourisboure et al.~\cite{dourisboure2007extraction}.

\section{Overview of methods for large scale optimization} \label{methods overview}

Nesterov \cite{Nesterov12} introduced constrained and unconstrained versions of an efficient method for solving convex  huge-scale optimization problems.
Followed by that paper,  appeared different versions of coordinate descent methods for large scale  convex optimization, see e.g., \cite{{Richtarik:14},Noecara:17,Nesterov14}.
In this paper,  we propose two variants of the random coordinate descent method to solve the D$k$S.
In this section we provide a brief overview of  algorithms from the literature, and describe those that are relevant to our work in more details.

The random coordinate descent method (RCDM) from \cite{Nesterov12}, is a method for solving unconstrained problems with convex objective.
The RCDM  performs in every iteration of the algorithm a random coordinate index selection  by using a  random counter.
The random counter generates  numbers   according to a distribution that is based on the coordinatewise Lipschitz constants.
The uniform coordinate descent method (UCDM) from \cite{Nesterov12} is developed for solving constrained problems with convex objective.
The method uses the uniform distribution to determine  random coordinates.
In the UCDM, each coordinate  update  is based on a solution of an optimization subproblem.
The optimization subproblem considers  constraints of the original problem,
and takes care that the new point is in the vicinity of the previous one.
The methods introduced in  \cite{Nesterov12} turn to be efficient for solving huge scale convex optimization problems.

There exist several extensions of the  random coordinate descent method  and uniform coordinate descent method from \cite{Nesterov12}.
For example, the random block coordinate descent method for linearly constrained optimization  by Necoara, Nesterov and Glineur \cite{Noecara:17}.
This method is introduced for solving problems with  a separable convex objective function and one linear constraint.
Richt\'arik and Tak\'a\v{c} \cite{Richtarik:14} extend results from \cite{Nesterov12} to composite optimization.
In particular, they introduce  randomized block-coordinate descent methods for minimizing composite functions.
Another recently developed method   for solving  large-scale optimization problems is a subgradient method by Nesterov  \cite{Nesterov14}.
The approach from \cite{Nesterov14} is suitable for  optimization problems with sparse subgradients.
The above mentioned methods are tested on large or huge scale convex problems such as the Google's PageRank problem, the PageRank problem,   image processing,
estimation in sensor networks or distributed control, $l_{1}$-regularized least squares problems.

On the other hand, there are very few results on solving large-scale nonconvex problems.
Patrascu and Necoara \cite{PatraNecora} derive  random coordinate descent algorithms for large scale structured nonconvex optimization problems,
and test them on sparse instances of the eigenvalue complementarity problem.

\medskip

Before we outline the UCDM  from  \cite{Nesterov12} and the 2-random coordinate descent algorithm from \cite{PatraNecora}, we introduce the notation.
Consider the space $\R^N$, and its decomposition on $n$ subspaces where $N=\sum_{i=1}^n n_i$.
We denote a block decomposition of $N\times N$ identity matrix by $I_N=(U_1, \ldots, U_n)\in R^{N\times N}$, where $U_i\in \R^{N\times n_i}$ ($i=1,\ldots,n$).
Thus, for $x=(x^{(1)}, \ldots, x^{(n)}) \in \R^N$ we have
\begin{equation}\label{decmposition}
x=\sum\limits_{i=1}^n U_i x^{(i)}
\end{equation}
where $x^{(i)} \in \R^{n_i}$  for $i=1,\ldots, n$. Note that similar notation is used in the related literature,  see e.g., \cite{Nesterov12}.

Let us now describe  the UCDM from \cite{Nesterov12}.
Consider a function  $f(x)$ that is  convex and differentiable on a closed convex set $Q\subseteq \R^N$. Assume that
the gradient of $f$ is coordinatewise Lipschitz continuous with constants $L_i$ ($i=1,\ldots, n$) where
\begin{equation} \label{Lj}
||  \nabla_{i} f(x+U_i h_i) -  \nabla_{i} f(x) ||\leq L_{i} || h_{i}|| \quad h_i\in \R^{n_i}, \quad i=1,\ldots, n, \quad x\in \R^N,
\end{equation}
$\nabla_i f(x)$ denotes the partial gradient of $f(x)$ in $x^{(i)}$, i.e.,
\[
\nabla_i f(x) = U_i^{\mathrm T} \nabla f(x) \in \R^{n_i}, \quad x\in \R^N,
\]
 and $\| \cdot \|$ denotes  the Euclidean norm.

Now,  the constrained optimization problem considered in \cite{Nesterov12}  is:
\[
\min_{x\in Q} f(x),
\]
where $Q= \bigotimes_{i=1}^n Q_i$ and the sets $Q_i \subseteq \R^{n_i}$ ($i=1,\ldots, n$) are closed and convex.
The  $i$th ($i=1,\ldots, n$) constrained coordinate update from \cite{Nesterov12} is:
\begin{equation} \label{upateV}
 V_{i}(x) = x + U_{i}^{\mathrm T}(u^{(i)}(x) - x^{(i)}),
\end{equation}
where and $u^{(i)}(x)$ is the solution of the following optimization subproblem:
\begin{equation} \label{upateU}
u^{(i)}(x) = {\rm arg} \min_{u^{(i)} \in Q_i}  \left [ \langle  \nabla_i f (x), u^{(i)} - x^{(i)} \rangle + \frac{L_{i}}{2} \| u^{(i)} - x^{(i)} \|^{2} \right ].
\end{equation}
Here $\langle  \cdot, \cdot \rangle$ denotes a vector product.
The uniform coordinate  descent method \cite{Nesterov12} chooses a random number $i$ from the discrete  uniform distribution and updates  $x^{(i)}$ in every iteration.
 In particular,  see Algorithm \ref{alg1}.

\begin{algorithm}%[t]
\caption{Algorithm UCDM \cite{Nesterov12}} \label{alg1}
\begin{algorithmic}
\REQUIRE A feasible initial solution $x_{0}$.
\STATE $k \leftarrow 0 $
\LOOP
\STATE Choose randomly $i_{k}$ by uniform distribution on $\{1,...,n\}$.
\STATE   Update $x_{k+1} = V_{i_{k}}(x_{k})$, by using \eqref{upateV} and \eqref{upateU}.
\STATE $k \leftarrow k+1 $
\ENDLOOP
\end{algorithmic}
\end{algorithm}

%\begin{algorithm}%[h!] %% tis one works with another packuage
%\caption{Algorithm UCDM \cite{Nesterov12}}
%%\SetAlgoLined
% %Start with a feasible initial solution $x_{0}$.
%\KwIn{A feasible initial solution $x_{0}$}
%%\For{$k \geq 0$}{
%\Repeat{}{
%Choose randomly $i_{k}$ by uniform distribution on $\{1,...,n\}$\\
%Update $x_{k+1} = V_{i_{k}}(x_{k})$, by using \eqref{upateV} and \eqref{upateU}.
%}{}
%%\Loop{Choose randomly $i_{k}$ by uniform distribution on $\{1,...,n\}$.}
%%\Endloop
%%\LOOP
%%\STATE Choose randomly
%%\ENDLOOP
%\end{algorithm}

The improvement in each step of the UCDM  is as follows:
\[
f(x) - f(V_{i}(x))  \geq \frac{L_{i}}{2} \| u^{(i)}(x) - x^{(i)} \|^{2}.
\]

\medskip
Patrascu and Necoara \cite{PatraNecora} introduce random coordinate descent algorithms for large-scale structured nonconvex optimization problems.
They consider unconstrained and linearly constrained problems with a nonconvex and composite objective function.
In particular, in \cite{PatraNecora} it is considered the following linearly constrained   optimization problem:
\[
\begin{array}{rl}
\min\limits_{x\in \R^N} & g(x)+ \l(x)\\[1ex]
{\rm s.t.} & a^{\mathrm T}x=b,
\end{array}
\]
where $a\in \R^N$ is a nonzero vector, $b\in \R$,  $g$ is a smooth function, and $\l$ is a convex, separable, nonsmooth function.
Further, the function $g$ has 2-block coordinate Lipshitz continuous gradient, i.e., there exist constants $L_{ij}>0$ such that
\[
||  \nabla_{ij} ~g(x+U_i h_i + U_j h_j) -  \nabla_{ij} ~g(x) ||\leq L_{ij} ||h_{ij}||
\]
for all $h_{ij}= [h_i^{\mathrm T}, h_j^{\mathrm T}|^{\mathrm T} \in \R^{n_i+n_j}$, $x\in \R^N$ and $i,j=1,\ldots, n$.
For given a feasible initial point $x_0$, that is $a^{\mathrm T}x_0=b$, the 2-random coordinate descent algorithm from \cite{PatraNecora}  is presented as Algorithm \ref{alg2}.\\

\begin{algorithm}%[t]
\caption{Algorithm 2-RCD \cite{PatraNecora}} \label{alg2}
\begin{algorithmic}
\REQUIRE A feasible initial solution $x_{0}$.
\STATE $k \leftarrow 0 $
\LOOP
\STATE  Choose randomly 2 block coordinates $(i_{k},j_k)$ with probability $p_{i_k,j_k}$.
\STATE  Update $x_{k+1} = x_k + U_{i_k}d_{i_k} + U_{j_k}d_{j_k}$.
\STATE $k \leftarrow k+1 $
\ENDLOOP
\end{algorithmic}
\end{algorithm}

Here, directions $d_{i_k,j_k}=[d_{i_k}^{\mathrm T},d_{j_k}^{\mathrm T}]^{\mathrm T}$
are obtained from the following optimization subproblem
\[
\begin{array}{rcll}
d_{i_k,j_k} &=& {\rm arg}\min\limits_{h_{i_k,j_k}} & g(x_k)+ \langle \nabla_{i_k,j_k} ~g(x_k), h_{i_k,j_k}  \rangle +
   \frac{L_{i_k,j_k}}{2} ||h_{i_k,j_k} ||^2 + \l (x_k+ h_{i_k,j_k}) \\[1.5ex]
 && {\rm s.t.} & a_{i_k}^{\mathrm T} h_{i_k} +  a_{j_k}^{\mathrm T} h_{j_k} = 0.
\end{array}
 \]
In \cite{PatraNecora}, the authors prove asymptotic convergence of the  sequence generated by 2-RCD to stationary points.
The authors  suggest updating only two block coordinates simultaneously in each iteration of the algorithm.
Thus,  the  2-RCD algorithm  updates only two coordinates simultaneously when $n=N$ i.e.,  for the scalar case.
In the following section we consider updating more than two coordinates in each iteration of our algorithms.
Our numerical results show that updating more than two (block) coordinates simultaneously  results with a very efficient algorithm.

\section{Two new algorithms for solving the D$k$S} \label{sect:onSolveDkS}

We present here two new algorithms for solving  large instances of the D$k$S.
The algorithms are tailored for solving the relaxation of the D$k$S, but can be adjusted for solving any nonconvex problem with a linear constraint.
Our algorithms can be seen as extensions of the UCDM and the 2-RCD algorithms, see Section \ref{methods overview}.
While our first algorithm converges  to a real feasible point,
our second algorithm provides an integer feasible point in most of the cases, see Section \ref{sect:numer}.
A version of the here presented second algorithm was studied in a master thesis by  van der Doef \cite{SusaneDoef}.

Let us first consider the following relaxation of \eqref{QP}:
\begin{equation}\label{lpQP}
\begin{array}{rl}
\max & x^{\mathrm T}Ax \\[1ex]
{\rm s.t.} & \sum\limits_{i=1}^{n} x_{i} = k \\[1.5ex]
& 0\leq x_{i} \leq 1, ~~ \forall i \in \{1,...,N \}.
\end{array}
\end{equation}

Recall that the constrained coordinate update in the  UCDM algorithm  considers one block coordinate,
while the coordinate update in the  2-RCD algorithm updates two block coordinates in each iteration.
To solve  the D$k$S we update  several coordinates simultaneously, in each iteration of our algorithms.

We set $f(x)= x^{\mathrm T}Ax$ and suppose that $n=N$, see \eqref{decmposition}.
In each step of our coordinate descent algorithms we update $q \geq 2$ coordinates.
Let $J_{i}$, $|J_i|=q$, be the set of random coordinates that are updated simultaneously in step $i$.
Then, our $q$-random coordinate constrained  update in $i$th  iteration is as follows:
\begin{equation}  \label{upateW}
W^i(x)_j = \left \{
\begin{array}{ll}
u^{i}_j(x) & {\rm if~} j\in J_{i} \\
x_j & {\rm otherwise}
\end{array} \qquad j=1,\ldots, n,
\right .
\end{equation}
where  $u^{i}(x)\in \R^q$  is the solution of a concave optimization subproblem. In particular,
\begin{equation}  \label{upateUnew1}
\begin{array}{rlrl}
u^{i}(x)&=&  {\rm arg}\max\limits_{u^{i}} &
 \sum\limits_{j\in J_i}   \nabla_j f (x) (u^i_j- x_j) - \sum\limits_{j\in J_i}  \frac{L_{j}}{2}  (u^i_j- x_j)^2   \\[2ex]
&&\text{s.t.} &   \sum\limits_{j\in J_i} u_j^{i} = k - \sum\limits_{j\notin J_i}  x_j \\[2ex]
&&& 0 \leq u^{i}_{j} \leq 1 \quad  \forall j \in J_{i}.
\end{array}
\end{equation}
Here,  $L_j$ is defined as in \eqref{Lj}.
Note that the  $q$-random coordinate constrained  update in $i$th  iteration can be also written as:
\[
W^i(x)=  x + \sum_{j\in J_i} U_j(u^i_j(x)-x_j).
\]
It is a well known result that  convex quadratic problems are polynomially solvable, see \cite{Kozlov}.
Kozlov et al.~\cite{Kozlov} reported $O(n^4L)$ algorithm for convex quadratic problems,
where $n$ is the number of variables and $L$ is the size of the problem.
Later papers present algorithms that have complexities of  $O(n^{3}L)$ arithmetic operations, see e.g., \cite{Goldfarb,Kojima}.

Now, we are ready to show our first algorithm.
For a fixed $q$ ($2\leq q \leq n$)  the $q$-random coordinate constrained algorithm $q$-RCC1  is presented as Algorithm \ref{alg3}.\\
\begin{algorithm}%[t]
\caption{Algorithm $q$-RCC$1$} \label{alg3}
\begin{algorithmic}
\REQUIRE A feasible initial solution $x_{0}$.
\STATE $k \leftarrow 0 $
\LOOP
\STATE  Determine $J_{k}$:  choose $q$ coordinates randomly by uniform distribution on $\{1,...,n\}$.
\STATE  Update $x_{k+1} = W^{k}(x_{k})$ by using \eqref{upateW} and \eqref{upateUnew1}.
\STATE $k \leftarrow k+1 $
\ENDLOOP
\end{algorithmic}
\end{algorithm}

Note that one can solve \eqref{upateUnew1} efficiently with a convex quadratic programming solver.

Since the optimization problem \eqref{lpQP} is nonconvex, the algorithm $q$-RCC$1$ can stuck in a local optimum.
Therefore, we also allow restarting of  the algorithm  from a new feasible starting point.
The algorithm $q$-RCC$1$ uses several stopping criteria.
For details on restarting and stopping criteria, see Section \ref{sect:numer}.
Numerical results show that  $q$-RCC$1$ converges to a local optimum of the  relaxation problem  \eqref{lpQP}.
However, we are interested in solving the integer programming problem \eqref{QP}.
Therefore, the subproblem of our next algorithm  considers only a linear approximation of the nonconvex objective  function.
In particular, our second algorithm  solves the following subproblem in order to find a $q$-random coordinate constrained update in $i$th iteration:
 \begin{equation}  \label{upateUnew2}
\begin{array}{rlrl}
u^{i}(x)&=&  {\rm arg}\max\limits_{u^{i}} &
   \sum\limits_{j\in J_i}   \nabla_j f (x) (u^i_j- x_j)   \\[2ex]
&&\text{s.t.} &   \sum\limits_{j\in J_i} u_j^{i} = k - \sum\limits_{j\notin J_i}  x_j \\[2ex]
&&& 0 \leq u^{i}_{j} \leq 1 \quad  \forall j \in J_{i}.
\end{array}
\end{equation}
Thus, to obtain  the $q$-random coordinate constrained update  from \eqref{upateUnew2}, we need to solve a linear programming  problem.

In 1979,  Khachiyan \cite{Khachiyan79} proved that linear programming is polynomially solvable.
Karamarker's  well known projective algorithm, see \cite{Karamark}, solves linear programs  in $O(n^4L)$ operations,
where $n$ is the number of variables in a standard-form problem with integer data of bit size $L$.
Many subsequent papers have reported $O(n^3L)$ algorithms for linear programming.
Anstreicher \cite{Anstreicher} shows that the complexity to solve linear programming problems  can be reduced to $O([n^3/\ln n]L)$. 

Our  $q$-random coordinate constrained  algorithm  $q$-RCC2 for a fixed $q$ ($2\leq q \leq n$) is given as Algorithm \ref{alg4}.\\

\begin{algorithm}%[t]
\caption{Algorithm $q$-RCC$2$} \label{alg4}
\begin{algorithmic}
\REQUIRE A feasible initial solution $x_{0}$.
\STATE $k \leftarrow 0 $
\LOOP
\STATE  Determine $J_{k}$: choose $q$ coordinates randomly by uniform distribution on $\{1,...,n\}$.
\STATE  $x_{k+1} = W^{k}(x_{k})$ by using \eqref{upateW} and \eqref{upateUnew2}.
\STATE $k \leftarrow k+1 $
\ENDLOOP
\end{algorithmic}
\end{algorithm}

Our numerical results show that the algorithm   $q$-RCC$2$ converges to an integer point for sufficiently  large $q$.
Once the algorithm finds an integer point, all points in  successive iterations are also integer.
Moreover, the successive integer vectors might be  in the vicinity  of the first found integer solution.
Therefore, the $q$-RCC$2$ algorithm could end up again in the first found integer solution.
To prevent cycling, we stop the $q$-RCC$2$ algorithm once the first integer solution is found,
or we restart the algorithm  from a new feasible  point.

\section{Numerical Results} \label{sect:numer}

In this section we present numerical results on solving the D$k$S problem by using our two algorithms, i.e., $q$-RCC$1$ and $q$-RCC$2$.
Numerical results are performed on an  Intel Xeon, E5-1620, 3.70 GHz with 32 GB memory.
To compute \eqref{upateUnew1} (resp.~\eqref{upateUnew2} we use Cplex 12.6 QP  (resp.~Cplex 12.6 LP) solver.

We test our algorithms on random graphs and several graphs from the literature.
In particular, we consider the following graphs:
\begin{enumerate}
\item{\sc The Erd\"os-R\'enyi  graph:}  Each edge in a graph is generated independently of other edges with probability $p\in (0,1]$.
For any given $p$, a graph formulated in the described way is known  as the Erd\"os-R\'enyi random graph $G_p(n)$.
The Erd\"os-R\'enyi  graph was introduced by Erd\"os and R\'enyi  in 1959, see  \cite{ErRe:59,{ErRe:60}}.

\item{\sc The Erd\"os-R\'enyi  graph with a planted subgraph:}  In the   Erd\"os-R\'enyi  random graph  $G_p(n)$ we plant a complete subgraph with $k$ vertices.
We denote the resulted graph by $P_p^k(n)$.
Random graphs with planted subgraphs are used also in Tsourakakis et al.~\cite{tsourakakis2013}. We compare our results with their heuristic results.

{
\item{\sc Instances for the D$k$S from the literature:}
\begin{itemize}
\item  We consider instances that are available from the following webpage:\\
{\url{http://cedric.cnam.fr/~lamberta/Library/k-cluster.html}}.
Those instances are used as test instances for the densest $k$-subgraph problem in
\cite{billionnet2005different,billionnet2009improving,{BillionWiegele17},{KrisMalRoup:16}}.
The sizes of instances  are   $n=40$, $80$, $100$, $120$, $140$, $160$, and densities $d=25, 50, 75 ~\%$.
For a given number of vertices $n$ and a density $d$ an unweighted graph is randomly generated.
In all above mentioned papers, the parameter $k$ has following values: $\frac{n}{2}$, $\frac{n}{4}$ and $\frac{3n}{4}$.
 Therefore,  we also use the same values for $k$ in our experiments.

\item  Brimberg et al.~\cite{brimberg2009variable} generated test instances for the heaviest $k$-subgraph problem.
One can download those instances from the following page \\
{\url{http://www.mi.sanu.ac.rs/~nenad/hsp/}}.
There are 177 instances in total, for details see \cite{brimberg2009variable}.
%%We provide a comparison of our algorithms with  heuristic approaches from \cite{brimberg2009variable} on those instances.

\end{itemize}
}

\item{\sc Real-world data:} We consider data from the following two different  sources.
\begin{itemize}
\item We test our algorithms on several graphs from 10th DIMACS Implementation Challenge - Graph Partitioning and Graph Clustering.
In particular, we consider Jazz graph and Email graph.
These two graphs are derived from two different networks and then symmetrized, as explained on the  DIMACS webpage:
\url{https://www.cc.gatech.edu/dimacs10/archive/clustering.shtml}\\

\item  We take graphs from the following webpage: \url{snap.stanford.edu}.
In particular, we use  undirected graphs  from the collaboration networks database.
A collaboration network represents scientific collaborations between authors of papers in a specific field.
A graph from the database is represented by an adjacency matrix whose element on position  $(i,j)$
 equals one if author $i$ co-authored a paper with author $j$.
The largest here considered graph from the snap.stanford.edu database has 23,133 vertices.

\end{itemize}
\end{enumerate}

We present below settings of our algorithms:
\begin{itemize}
\item We tested two different types of initial feasible  points.
The first type of an initial point is a random point in which all coordinates have values between zero and one, and all sum up to $k$.
The second type of an initial point is the vector whose all coordinates equal to $k/n$.
Our numerical results show that the quality of a solution computed by our algorithms does not depend on a starting point.
Therefore, our algorithms start with a randomly generated feasible point unless indicated differently.
Namely, it is costly to use the first type of the starting point when $n>2^{13}$,
and therefore we use the second type of initial point for larger $n$.

\item We implement several stopping criteria.
Both algorithms $q$-RCC$1$ and $q$-RCC$2$ stop after a pre-specified number of iterations  is reached,
unless any other stoping criteria is satisfied. We list the remaining stoping criteria below.

\begin{itemize}

 \item Stop when the first integer solution is found. Our tests show that after the first integer solution is found,
 $q$-RCC$2$ might cycle i.e., end up in the same integer point after a certain number of iterations.
 This happens since the algorithm computes integer points the vicinity of the first found integer solution.
 Our tests show that it is better to stop the algorithm when the first integer solution is found,
 instead of letting it run till eventually cycling  appears and then stop.
 Namely, the latter requires more computational effort, but does not necessarily  result in a significant improvement of the objective value.
 On the other hand,  $q$-RCC$1$ does not converge to an integer solution, in general.
 However, for large $q$ the algorithm  $q$-RCC$1$  might also provide an integer solution.
 In the case that $q$-RCC$1$ finds an integer point, we stop the algorithm.

\item Stop $q$-RCC$1$ if the difference in two consecutive objective values is less than a pre-specified tolerance.  We use here $\epsilon = 1e-7$ as the tolerance.
  This criteria is not implemented in $q$-RCC$2$ since the algorithm tends to find faster  an integer value than to satisfy this criterion.

\end{itemize}

\item Restarting of the algorithms.
We sometimes perform restarting of the algorithms $q$-RCC$1$ and $q$-RCC$2$  for a given number of times and
after one of the stopping criteria from above is reached.
In each new run, we restart the algorithm  by using one of the previously described starting points.
It is going to be clear from the context if we performed restarting of the algorithm.

\end{itemize}

It might happen that in  an iteration  of our algorithm, the objective value decreases and then in the next iterations keeps improving.
This happens since the objective is nonconvex. Extensive test shows that there is no harm in accepting non-improving moves,
since the algorithms recover fast.
We tested  our algorithms also  when only improving moves are accepted, and concluded that there is no benefit of doing this.

\medskip

Let us now present computational  results.
 We test our two algorithms on various instances, and present lower bounds for the problem
 \eqref{lpQP} obtained from the limit point returned by the algorithms. 

\medskip
\noindent
{\bf Tests on the Erd\"os-R\'enyi  graphs}. \\
We first show performance of our algorithms on $G_{0.5}(2^{10})$ for different number of simultaneous updates $q$ and different number of iterations.

Table \ref{tabRand1} presents bounds computed  by  the algorithm  $q$-RCC$1$  for one graph only with 1,024 vertices and   $k=30$.
Here, we do not restart  $q$-RCC$1$.
The initial point in all runs have coordinates $k/n$ with objective value  622.54.
The table reads as follows. In the first row we specify $q$.
Rows indicated by bnd.~provide bounds that are computed in  seconds, given  in the first row below that one.
Finally,  rows indicated by iter.~specify the number of iterations needed to compute bounds listed in the two rows above that one.

All computations in Table \ref{tabRand1} terminated after  the algorithm performed a pre-specified number of iterations.
The results in Table \ref{tabRand1}  show that  the quality of bounds improve and corresponding computational times increase along with the number of iterations.
The table also shows that for $q=2$ there is a small improvement in the bound even after 10,000 iterations.
Note also that  for large number of updates i.e., $q=750$ there is no significant  improvement in the bound value when the number of iterations increases.
Table \ref{tabRand1} also indicates that a good strategy for computing bounds
is to take $q$ that is between 10$\%$ and 20$\%$ of the number of vertices in the graph.
\begin{table}[t]
\begin{center}
\begin{tabular}{c|ccccccc}
\hline
 %    & \multicolumn{8}{c|} {$q$} \\
 $q$     & 2        & 50     & 100       & 200     & 500     & 750 \\[1ex]
\hline
bnd. &  623.12  & 653.15 &  686.54   & 738.260 &  808.38 &  816.86 \\
time &  1.01    & 1.02   &    1.28     & 2.20   &   4.77 &  7.51    \\
iter.&  500     &  500   &    500    &  500    &    500  &  500 \\[3ex]

bnd. & 623.67   & 682.96 &    741.23  &  795.44 &  825.81 &  837.91 \\
time & 1.95     &  2.10  &    2.56    &   9.56 &   9.56 &  15.18  \\
iter.&  1000    &  1000  &    1000    &   1000  &    1000 &  1000 \\[3ex]

bnd. & 628.12  &  799.33 &   824.27  &   831.99 &   836.96 &   833.33 \\
time &  9.54  &   10.33  &   13.32   &   23.66  &   49.27 &  77.55    \\
iter.&  5000    & 5000   & 5000  & 5000  & 5000  & 5000  \\ [3ex]

bnd. & 634.59 & 819.93 &  837.06   & 842.81  & 837.99  &  840.83  \\
time & 16.55  &  21.02  & 27.35   &  48.50 &  102.74 &    154.94 \\
iter.& 10000  & 10000   & 10000    &   10000  &   10000 &    10000 \\ \hline
\end{tabular}
\caption{$q$-RCC$1$ for $G_{0.5}(1024)$: bounds, running times (s) and iterations.}  \label{tabRand1}
\end{center}
\end{table}

Table \ref{tabRand2}  presents bounds computed by  $q$-RCC$2$ for the same graph used in  Table \ref{tabRand1}.
Since the algorithm  $q$-RCC$2$ with $q>2$  terminates in most of the cases due to the  stopping criteria
``the first integer solution is found",
we present results obtained by averaging  20 bounds computed after  20 times restarting the algorithm with the same starting point.
We present average of 20 runs for each $q$.
The initial point in all runs is the vector with  coordinates $k/n$.
Table \ref{tabRand2} shows that the average of 20 bounds is the best for $100$ simultaneous updates.
Note that for $q\geq 200$ the computational time significantly drops, but the quality of  bounds deteriorate.
Among all computed bounds  the best integer value is 840.
That  value is obtained for $q=50$ and for $q=100$.

\begin{table}[t]
\begin{center}
\begin{tabular}{c|ccccccc}
\hline
 %    & \multicolumn{8}{c|} {$q$} \\
 $q$     & 2        & 50     & 100       & 200     & 500     & 750 \\[1ex]
\hline
bnd. & 636.15 & 825.23 &  828.66 &     824.20  & 790.70  & 731.30   \\
time & 0.643  & 1.12   &  1.07   &       0.64  &   0.26  &   0.20   \\
iter.& 500    & 492.05 &  255.70 &     70.90   &  13.9   &  6.85     \\ \hline
\end{tabular}
\caption{$q$-RCC$2$ for $G_{0.5}(1024)$: bounds, running times (s) and iterations.}  \label{tabRand2}
\end{center}
\end{table}

\medskip\medskip
\noindent
{\bf Tests on graphs with  planted subgraphs}. \\
We plant complete subgraphs in random graphs  because we know the optimal value of the problem. This enables us to evaluate the performance of our algorithms.
Note that  heuristic approaches \cite{brimberg2009variable,tsourakakis2013} report results for graphs with up to 3,000 vertices.

Table \ref{tablPlanted1} summarizes outcomes of our two algorithms on graphs with planted subgraphs and 4,096 vertices.
In particular, we consider the  Erd\"os-R\'enyi   graphs  $G_{0.3}(2^{12})$ whose planted complete subgraphs have  $100$ vertices.
Note that the optimal value of the densest 100-subgraph problem on the described graph  is 9,900 with high probability.
We run each algorithm with different $q$ on 30 different $P_{0.3}^{100}(4096)$  graphs.
In particular, we run 1,000  iterations of $q$-RCC$1$  and  1,000  iterations of $q$-RCC$2$  for each  $q\in \{2, 400, 800, 2000\}$ and each graph.
In the row denoted by $q$-RCC$1$ (resp.~$q$-RCC$2$) we list the best obtained bound among 30 values for the given $q$, as well as the computational time in seconds needed to compute that bound.

An interesting result is that the algorithm $q$-RCC$1$ with  $q=2000$  computes the value $9,899.99$  for 28 different graphs.
Coordinates of the solution vectors in those 28 cases differ at most for {$1e-5$} from the value 0 or 1.
If we let run  $q$-RCC$1$  with  $q=2,000$ for 10,000 iterations  the best obtained result is $9,899.99996$.
Here, values of coordinates in the solution vector are within an error of {$1e-6$} from  0 or 1. It takes 1,215 seconds to perform 10,000 iterations.

The algorithm $q$-RCC$2$ finds the  value  9,900 in  11, 8 and 13 cases  for $q=400$,  $q=800$, $q=2,000$, respectively.
Table \ref{tablPlanted1} reports the shortest computational time required to compute  9,900 by   $q$-RCC$2$ among all computations.
The longest time needed to obtain 9,900   by $q$-RCC$2$ and $q=800$ (resp.~$q=2,000$) is 35.03 s  (resp.~72.64 s).
Finally,  $q$-RCC$2$ computes the weakest bound for $q=2,000$.
The results in Table \ref{tablPlanted1} show that the algorithm $q$-RCC$2$  performs  better than $q$-RCC$1$   for all $q$.
However, $q$-RCC$2$  can stop fast in a weak bound.
On the other hand,  $q$-RCC$1$ improves slowly and steadily.
\begin{table}[t]
\begin{center}
\begin{tabular}{c|cccc}
\hline
 $q$          &  2                  & 400                 & 800               & 2000    \\[1ex]\hline
$q$-RCC$1$  & 6185.23 (2.08)    & 6908.02 (24.63)    & 9872.83 (46.30) & 9899.99 (121.53)  \\[1ex] %\hline <--doublechecked times
$q$-RCC$2$  & 6240.01 (0.43)     & 9900 (5.83)        &  9900 (2.42)     &   9900 (1.26) \\ \hline
%%$q$-RCC$2$  & 6225.99 (0.48s)     & 7814 (10.19s)       &  7770 (2.50s)    &   7648 (1.03s) \\ \hline
%%          &                     & 10 times            & 8 times          & 9 times
\end{tabular}
\caption{Bounds  and running times (s) for  $P_{0.3}^{100}(4096)$.}  \label{tablPlanted1}
\end{center}
\end{table}

Let us now consider a similar experiment as the previous one, see Table \ref{tablPlanted2}.
In particular, in the  Erd\"os-R\'enyi   graph  $G_{0.3}(2^{12})$ we plant
a complete subgraph with $800$ vertices, which results in $P_{0.3}^{800}(4096)$.
Note that the optimal value of the densest 800-subgraph problem on the described graph is 639,200  with high probability.
We run 1,000  iterations of each of the algorithms for 30 different graphs and for $q=2, 400, 800, 2000$.
The algorithm $q$-RCC$2$ finds the value 639,200  in 28, 28 and 27 cases  for $q=400$,  $q=800$, $q=2000$, respectively.
The results in Table \ref{tablPlanted2}  indicate that the algorithm $q$-RCC$2$ finds faster and more frequently cliques with $800$ vertices than cliques with 100 vertices.

Finally, in a similar experiment  with 2000-planted subgraph problem,  the  optimal
 value  is computed by  $q$-RCC$2$ in  23, 28 and 27 cases  for $q=400$,  $q=800$, $q=2000$, respectively.
\begin{table}[t]
\begin{center}
\begin{tabular}{c|cccc}
\hline
  $q$         &  2                  & 400                 & 800               & 2000    \\[1ex]\hline
$q$-RCC$1$  & 215,748.73   (1.97)    & 639,199.99 (25.05)    &  639,199.99 (46.48) &  639,199.99  (121.07)  \\[1ex]
$q$-RCC$2$ & 228,805.97 (0.46)    &  639,200 (2.61)      &  639,200 (2.33)    &  639,200 (0.91)  \\ \hline
 %number    &  1            &   28  & 28 &  27  \\
\end{tabular}
\caption{Bounds and running times (s) for  $P_{0.3}^{800}(4096)$.}  \label{tablPlanted2}
\end{center}
\end{table}

We did also extensive tests on $G_{0.2}(2^{13})$  with planted cliques on $500$ vertices, i.e., $P_{0.3}^{500}(2^{13})$.
It is interesting to note that for those graphs and 800 simultaneous updates, $q$-RCC$2$ always finds the planted subgraph between 10 and 55 seconds.\\

Finally, we experiment with random graphs as in \cite{tsourakakis2013}.
We plant a complete graph with 30 vertices in $G_{p}(3000)$ with $p\in \{0.008, 0.1, 0.5 \}$.
The  algorithm $q$-RCC$2$ with $q=100$ finds the clique within 6 seconds in  $P_{0.008}^{30}(3000)$.
The  algorithm $q$-RCC$2$ with $q=150$ finds the clique within 50 seconds in $P_{0.01}^{30}(3000)$.
There are no computational times reported in \cite{tsourakakis2013}.
However, in \cite{tsourakakis2013}, the authors report that all considered algorithms find the clique in a graph $G_{0.008}(3000)$,
and only one algorithm can find the clique in  $G_{0.01}(3000)$.
On the other  hand, no algorithms from  \cite{tsourakakis2013} could find the clique in  $G_{0.5}(3000)$.
We also couldn't find the clique in $G_{0.5}(3000)$, even after several restarting of the algorithm.

\medskip\medskip
\noindent {
{\bf Tests on instances from the literature.} \\
We consider instances from  {\url{http://cedric.cnam.fr/~lamberta/Library/k-cluster.html}}.
Those instances are also used as test instances for the  D$k$S  in
\cite{billionnet2005different,{billionnet2009improving},{BillionWiegele17},{KrisMalRoup:16}}, see also
\url{http://www-lipn.univ-paris13.fr/BiqCrunch/results}.
We summarize  the outcome of our computational experiments below.

For each instance with $n=40$, any given density  and any  $k=10,20,30$,  our algorithm $q$-RCC$2$ finds an optimal solution within 0.1 s.
In  $q$-RCC$1$ we implement additional stopping criteria, that is to stop the algorithm when the objective value differs from
the optimal objective value for less than 0.0001. The  algorithm  $q$-RCC$1$  provides such bounds within  2 s.
We allow 1,000 iterations per round in both algorithms. Solutions are mostly found in the first round of the algorithms.
We test both  algorithms for $q=4$ and $q=8$ and notice similar performance of the algorithms for both values of $q$.

For each instance  with $n=80$, any given density  and any $k=20,40,60$,  the  algorithm $16$-RCC$2$ finds an optimal solution within 0.2 s,
while  $8$-RCC$2$  needs at most $0.4$ s.
The algorithm $16$-RCC$1$  performs better than $8$-RCC$1$ and requires at most 13 seconds to obtain a bound that is close to the optimal solution.
In most of the  cases, $16$-RCC$1$ finds an optimal solution  in less than 4 s.
Here, we use the same additional stopping criteria as for instances with 40 vertices. We allow 2,000 iterations per round.

For  each instance   with $n=100$, any given densities,  and  any $k=25,50,75$ the  algorithm $10$-RCC$2$ finds an optimal solution in less than 1 s.
There are several instances for which the algorithm runs up to 3 seconds.
We set 3,000 for the maximal number of iterations in one round.
This enables $20$-RCC$1$ to converge to an optimal solution of a given instance  in at most 16 seconds.

For instances with $n=120,140, 160$ we tested only the algorithm $q$-RCC$2$.
We take for $q$ the value that is equal to  $20\%$ of the number of vertices in the given instance, and allow 3,000 iterations per round.
For most of the instances   $q$-RCC$2$ finds  optimal solutions within 4 seconds.
For the instances  {\tt kcluster160-050-40-1.dat}, {\tt kcluster160-050-40-5.dat}, {\tt kcluster160-075-40-2.dat} and {\tt kcluster160-075-80-4.dat}
we needed to change the value of $q$ in order to find optimal solutions.
In particular,  we set $q$ to be $15\%$ of the number of vertices and found optimal solutions within  25 seconds. \\

\medskip\medskip

Brimberg et al.~\cite{brimberg2009variable} provide extensive computational experiments on solving the heaviest $k$-subgraph problem by using several heuristic approaches.
In particular, they compare  performances of the following heuristics: two greedy constructive heuristics (drop and add),
 two versions of variable neighbourhood search  (VNS)  heuristics (basic VNS and skewed VNS),
two tabu search heuristics (TS1 and TS2) and two multi-start heuristics (MLS1 and MLS2). The results in \cite{brimberg2009variable} show
that VNS heuristic preforms the best over other heuristics. On the other hand TS1 has  the worst performance among tested approaches.

%%In Table \ref{heaviest} we also present time required to solve instances by $q$-RCC$2$.
%% Since running  time of VNS and TS1 is taken from \cite{brimberg2009variable},
%% the reader should take into account that the machine used in those computations is weaker than ours.

Here, we test  the  algorithm $q$-RCC$2$ on the same set of instances as in \cite{brimberg2009variable}.
We compare our results with the VNS and TS1 heuristics that use random initial starting points, see Table 3 in \cite{brimberg2009variable}.
In  Table \ref{heaviest} we report the average $\%$ deviation
\[
\% ~{\rm deviation} = \frac{ {\rm best~value} - {\rm rcc2}   }{ {\rm best~value} } \cdot 100,
\]
where  `best value' denotes the best known solution reported in \cite{brimberg2009variable}, and `rcc2' denotes our bound.
We also report average running time obtained by $q$-RCC$2$, see the last column in Table \ref{heaviest}.
To solve instances  we set $q=100$ and alow restarting the algorithm 100 times.
For instance with 1,000 nodes we set 7,000 for the maximum number of iterations per round,  while for instances with  3,000 nodes we set 10,000 iterations per round.

The results in Table \ref{heaviest} show that our algorithm is performing  better than TS1 and worse than  VNS.
Note that  the average $\%$ deviation of our algorithm is within $2\%$.
The algorithms from \cite{brimberg2009variable} are specialized for solving the heaviest $k$-subgraph problem, while we use the best settings for the D$k$S.

%
%\begin{table}[t] \label{heaviest}
%\begin{center}
%\begin{tabular}{ccccccccc}
%\hline
% & & & \multicolumn{3}{c}{ } & \multicolumn{3}{c}{ } \\
%
%type    & $n$     & $k$   &  $\%$ deviation & & &time &      \\
% &&& $q$-RCC$2$ & VNS & TS1  &     $q$-RCC$2$ & VNS & TS1 \\[1ex]\hline
%%\multicolumn{4}{|c|}{ } \\
%
%I sparse  &  1000 & 300 &  0.90  & 0.15 & 1.64 & 113.05 & 72.20 & 3.13\\[1ex]
%I sparse  &  1000 & 400 &  0.55  & 0.10 & 0.99  & 103.64 & 75.78 & 3.86  \\[1ex]
%I sparse  &  1000 & 500 &  0.24   & 0.03 & 0.55 &  89.90 & 26.14 & 3.99 \\[3ex]
%
%I dense  & 1000 & 300 &   0.24 & 0.04 & 0.49 & 112.91 & 225.47 & 85.90 \\[1ex]
%I dense  & 1000 & 400 &  0.08 & 0.03 &0.40 &  90.15 &  52.67 & 3.98\\[1ex]
%I dense  & 1000 & 500 & 0.03 & 0.00 & 0.20 &  131.65 & 185.41 & 4.30 \\[3ex]
%
%II sparse  &  3000 &  900 & 1.09 & 0.05 & 1.52 &  161.18 &566.53 & 191.45 \\[1ex]
%II sparse  &  3000 & 1200 & 0.59 & 0.03 & 1.02 & 238.37 & 349.59 & 217.21 \\[1ex]
%II sparse  &  3000 & 1500 & 0.30 & 0.00 & 0.57 &  197.37 & 376.47 &231.80 \\[3ex]
%
%III sparse  &  1000 & 300 &   & 0.07  & 5.19  &  & 151.48  & 2.75  \\[1ex]
%III sparse  &  1000 & 400 &   & 0.04  & 2.36  &  & 118.98  & 3.76  \\[1ex]
%III sparse  &  1000 & 500 &   & 0.02  & 1.40  &  & 136.03  & 3.19 \\[3ex]
%
%\hline
%\end{tabular}
%\caption{Summary results for all three types of the heaviest $k$-subgraph problem.}  \label{tablPlanted2}
%\end{center}
%\end{table}

\begin{table}[t] \label{heaviest}
\begin{center}
\begin{tabular}{ccccccccc}
\hline
 & & & \multicolumn{3}{c}{ } & \multicolumn{1}{c}{ } \\

type    & $n$     & $k$   &  $\%$ deviation & & &time       \\
 &&& $q$-RCC$2$ & VNS & TS1  &     $q$-RCC$2$\\[1ex]\hline
%\multicolumn{4}{|c|}{ } \\

I sparse  &  1000 & 300 &  0.90  & 0.15 & 1.64 & 113.05 \\[1ex]
I sparse  &  1000 & 400 &  0.55  & 0.10 & 0.99  & 103.64  \\[1ex]
I sparse  &  1000 & 500 &  0.24   & 0.03 & 0.55 &  89.90  \\[3ex]

I dense  & 1000 & 300 &   0.24 & 0.04 & 0.49 & 112.91  \\[1ex]
I dense  & 1000 & 400 &  0.08 & 0.03 &0.40 &  90.15 \\[1ex]
I dense  & 1000 & 500 & 0.03 & 0.00 & 0.20 &  131.65  \\[3ex]

II sparse  &  3000 &  900 & 1.09 & 0.05 & 1.52 &  161.18  \\[1ex]
II sparse  &  3000 & 1200 & 0.59 & 0.03 & 1.02 & 238.37  \\[1ex]
II sparse  &  3000 & 1500 & 0.30 & 0.00 & 0.57 &  197.37  \\[3ex]

III sparse  &  1000 & 300 & 1.95  & 0.07  & 5.19  & 140.47  \\[1ex]
III sparse  &  1000 & 400 & 1.25  & 0.04  & 2.36  & 104.16   \\[1ex]
III sparse  &  1000 & 500 & 0.77  & 0.02  & 1.40  & 157.43  \\[3ex]

\hline
\end{tabular}
\caption{Summary results for all three types of the heaviest $k$-subgraph problem.}  \label{tablPlanted2}
\end{center}
\end{table}

}

\medskip\medskip
\noindent
{\bf Tests on real-world graphs.}\\
{\sc Jazz graph} represents jazz  musicians network related to $n=198$   musicians, see \cite{Gleiser}.
There are $m=2,742$ edges in the graph, which represent the network of jazz musicians.
It is known  that this graph contains a  clique with 30 vertices, see  \cite{tsourakakis2013}.
The algorithm $q$-RCC$2$ with $q=30$ finds the clique after  2 times restarting  the algorithm, which takes in total  0.07 seconds.
If we use $q=10$, then the algorithm finds the optimal clique after  10 times restarting  the algorithm.
For $q=2$ the algorithm $q$-RCC$2$ fails to find a clique even after restarting the algorithm 100 times. \\

{\sc Email graph} represents email network of $n=1,133$ members of the Univeristy Rovira i Virgili (Tarragona), see \cite{Guimera}.
There are $m=5,451$ edges in the graph.
From \cite{tsourakakis2013} we know that email graph  has a clique with 12 vertices.
Our algorithm $q$-RCC$2$  with $q=40$ finds the clique in 6 seconds (!). \\

Our final set of experiments consider graphs from the collaboration networks database.\\
{\sc CA-GrQc} collaboration network from \cite{Leskovec} covers scientific collaborations between authors of papers that are submitted to General Relativity and Quantum Cosmology category.
The data covers papers in the period of 124 months i.e.,  from January 1993 to April 2003. The adjacency matrix of the graph has 5,242 vertices and 14,496 edges.
We are not aware of an optimal value for the densest $k$-subgraph problem on {\sc CA-GrQc}.
Therefore, we present our results for different $k$, see Table \ref{tableCaGrQc1}.
In the row denoted by {\sc CA-GrQc} we list for each $k$ the best computed objective value and the corresponding computational time in seconds.
All results are obtained using the same settings: 200 simultaneous updates and 3,000 iterations per round.
Note that for $k=10,20,30,40$ we find the optimal cliques.
In all those cases we needed to restart the algorithm at most 4 times.
For $k=50$ we could not find a clique, and the best solution found is equal to 2,146.

Further we provide similar experiments for {\sc CA-HepTh} collaboration network, see  \cite{Leskovec}.
This network covers scientific collaborations between authors of papers submitted to High Energy Physics - Theory category.
The adjacency matrix of this graph is of order 9,877.
There are 	25,998 edges in this network. We set $q=1,000$ and look for the densest $k$-subgraph in the graph.
Again, we are not aware of an optimal value for the densest $k$-subgraph problem on {\sc CA-GrQc}.
Our computational results are given in Table \ref{tableCaGrQc1}. We find cliques for $k=10,20,30$.

{\sc CA-HepPh} collaboration network  considers scientific collaborations between authors of papers that are submitted to High Energy Physics - Phenomenology category, see  \cite{Leskovec}.
The data covers papers in the period of 124 months, i.e., from January 1993 to April 2003.   The resulted graph has 12,008 vertices and 118,521 edges.
The results for  the densest $k$-subgraph problem on {\sc CA-HepPh} for $k=10,20,40,50$ are given in  Table \ref{tableCaGrQc1}.
It is remarkable  that we can found cliques with 10, 20, 30, 40 and 50 vertices in short time. We use here $q=1,000$.

{\sc CA-AstroPh} collaboration network  covers scientific collaborations between authors of papers submitted to Astro Physics category, see  \cite{Leskovec}.
The data covers papers in the period from January 1993 to April 2003.
The adjacency matrix of the graph has 18,772 rows and 118,521 edges.  We take here $q=1,000$.
Our computational results show that for larger  $q$ the computations are too  expensive.
We find cliques for each $k$ in less than 13 minutes.
\begin{table}[t]
\begin{center}
\begin{tabular}{c|ccccc}
\hline
$k$          &  10  & 20 & 30 & 40 & 50   \\[1ex]\hline
{\sc CA-GrQc}    &  90 (0.3)  & 380 (0.4)   &  870 (0.6)   & 1560 (5.9) & 2146 (9.1)     \\
 {\sc CA-HepTh}  & 90 (35.4) &  380 (48.2) &  870 (152.5) & 1048 (176.7) &  1166 (396.2) \\
  {\sc CA-HepPh}  & 90 (212.0) & 380 (370.9)  & 870 (410.1)  &  1560 (74.8) & 2450 (38.6) \\
 {\sc CA-AstroPh} & 90 (369.5) & 380 (180.6) & 870 (529.5) & 1560 (447.5) & 2450 (748.4)  \\
 %%%{\sc CA-AstroPh} &  90 (159.7s) &  380 (3584.9s) &          \\
 \hline
 \end{tabular}
\caption{Results  obtained by  $q$-RCC$2$ for different $k$.}  \label{tableCaGrQc1}
\end{center}
\end{table}

Finally, we consider {\sc CA-CondMat}  collaboration network. This collaboration network considers
scientific collaborations between 23,133 authors whose  papers are submitted to Condense Matter category.
The resulted adjacency matrix has 	93,497 edges.
The data cover papers in the period from January 1993 to April 2003.
In {\sc CA-CondMat}  we found a clique with 10 vertices in 159  seconds,
and a clique with 25 vertices in 211 seconds. To find densest subgraphs we set  $q=1,300$.

\section{Conclusion}
There are many studies on  random coordinate descent algorithms for convex problems, but a very few results on solving nonconvex large scale problems.
In this paper we present two algorithms for solving large scale nonconvex problems with one linear constraint.
We exploit our algorithms to solve large scale instances of the densest $k$-subgraph problem.

The main difference between our algorithms and those in the literature is that we allow updating more than two random coordinates
simultaneously in each iteration of the algorithms.
Our numerical results demonstrate significant improvement in bounds for larger than two simultaneous updates of the algorithms, see Table \ref{tabRand1}--Table \ref{tablPlanted2}.
The $q$-RCC$2$ algorithm performs better than $q$-RCC$1$, and for an appropriate $q$ it converges to an integer solution of the problem.
Note that the $q$-RCC$2$ algorithm considers a linear approximation of the nonconvex objective function.

Our numerical results verify the efficiency of the here introduced approach.
For instance, we are able to find  densest $k$-subgraphs in real world graphs with up to  23,133 vertices in a few minutes.
  Our numerical results in Section \ref{sect:numer} can be used as a benchmark for solving the densest $k$-subgraph problem on large graphs.
   \\

\noindent
{\bf Acknowledgements.}
The author would like to thank Pavel Dvurechensky for useful discussions on the   random coordinate descent approaches.
The author would also like to thank two anonymous referees for suggestions that led to an improvement of this paper.

\end{document}